\documentclass{article}
\usepackage{a4}
\usepackage[latin1]{inputenc}
\usepackage[T1]{fontenc}
\usepackage{amsmath}
\usepackage{amssymb}
\usepackage{amsthm}
\usepackage{lmodern}
\usepackage[ps,dvips,arrow,matrix,tips]{xy}
\SelectTips{cm}{10}

\renewcommand{\phi}{\varphi}
\newtheorem{thm}{Theorem}[section]
\newtheorem{prop}[thm]{Proposition}
\newtheorem{lemma}[thm]{Lemma}
\newtheorem{conj}[thm]{Conjecture}

\newtheorem{defn}[thm]{Definition}
\newenvironment{pf}[1][]{\noindent{\emph{Proof}#1.} }{\hfill $\square$}

\DeclareMathOperator{\Spec}{Spec}
\renewcommand{\P}{{\mathbf{P}}}
\newcommand{\A}{{\mathbf{A}}}
\newcommand{\F}[1]{\mathbf{F}_{#1}}
\newcommand{\Falg}{\mathbf{F}}
\newcommand{\Fq}{{\F{q}}}
\newcommand{\N}{{\mathbf{N}}}

\newcommand{\Gal}{\mathrm{Gal}}
\newcommand{\tuple}[2]{#1, \mskip2.5mu \ldots \mskip-1mu, \mskip2.5mu #2}
\newcommand{\symm}[1]{\mathfrak{S}(#1)}
\newcommand{\abs}[1]{\left| \mskip1mu #1 \right|}

\title{Representing an element in $\Fq[t]$ as the \\sum of two irreducibles}
\author{Andreas O. Bender}

\begin{document}
\maketitle

\begin{abstract}
A monic polynomial in $\Fq[t]$ of degree $n$ over a finite field $\Fq$ of odd characteristic can be written as the sum of two irreducible monic elements in $\F{q}[t]$ of degrees $n$ and $n-1$ if $q$ is larger than a bound depending only on $n$. The main tool is a sufficient condition for simultaneous primality of two polynomials in one variable $x$ with coefficients in $\Fq[t]$. 
\end{abstract}
{\bf MCS: 11T55, 11R09}

\smallskip
\section{Introduction}

Which monic elements $F$ in a polynomial ring $\Fq[t]$ can be written as the sum of two monic irreducibles of unequal degrees not larger than the degree of $F$? In a nutshell, we prove that such a representation is possible in odd characteristic if $q$ is sufficiently large relative to the degree of the polynomial $F$. 

\begin{thm}\label{cor}
Let $\Fq$ be a finite field of odd characteristic and cardinality $q$ and let $F$ be a monic polynomial in $\Fq[t]$ whose degree $n$ is at least $2$. 

Then if 
$q> 3 n^{4}16^{n^4} (n+1)^{8n^4+1}$, 
there exist irreducible monic polynomials $F_1$ and $F_2$ in $\F{q}[t]$ with $\deg (F_1)=\deg(F)-1$ and $\deg (F_2)=\deg(F)$ such that
$$F=F_1 + F_2.$$
\end{thm}

The statement of the theorem above with a brief summary of the proof appeared already as theorem $2$ in~\cite{CRnote}, but note that an incorrect value of the lower bound on $q$ was given there. 

A significant part of the proof of theorem~\ref{cor} consists of an adaptation of the proof of the main theorem in~\cite{b-w}. 

\smallskip
The results above are motivated by a letter sent to Leonhard Euler and dated 7th June 1742~\cite{goldbach}, in which Christian Goldbach conjectured that every integer greater than $5$ is equal to the sum of three primes. In his answer~\cite{euler}, Euler noted that this conjecture is equivalent to the claim that every even integer greater than $2$ can be written as the sum of two primes. 

As for partial results towards the Goldbach conjecture, Vinogradov~(\cite{vin} or~\cite[chpt. 26]{davenport}) proved that every sufficiently large odd integer is the sum of three primes. Furthermore, Chen~(\cite{chen} or~\cite[chpt. 10]{nathanson}) showed that every sufficiently large even integer is the sum of a prime and a product of at most two primes. 

\smallskip
Goldbach's conjecture can also be formulated in the ring $\Fq[t]$ rather than in the ring of integers. We note that the absolute value of $F\in \Fq[t]$ is defined as the integer $\abs{F}=q^{\deg(F)}$. 
Furthermore, the analogon to an even integer is a polynomial in $\Fq[t]$ which is divisible by an irreducible polynomial of absolute value equal to $2$. The distinction between even and odd polynomials therefore exists only in $\F2[t]$. 

The function field version of the Goldbach conjecture can then be stated as follows. 

\begin{conj}[Conjecture 1.20 in~\cite{eff-hayes-book}]\label{polyngoldbach}
Let $\Fq$ be a finite field of characteristic $p$ and cardinality $q$ and let $F$ be a monic polynomial in $\Fq[t]$ of degree $n$ at least $2$ which is even if $q=2$ and not of the form $t^2+t+a$ if~$p=2$. 

Then there exist irreducible monic polynomials $F_1$ and $F_2$ in $\Fq[t]$ with $\deg (F_1)<\deg(F)$ and $\deg (F_2)=\deg(F)$ such that
$$F=F_1 + F_2.$$
\end{conj}

Car~\cite{car1} gave an upper bound in terms of $n$ and $q$ for the number of polynomials $F\in \Fq[t]$ of degree at most $n$ which are not the sum of two irreducibles whose degrees are at most equal to the degree of $F$. She also proved a function field version of Chen's theorem~\cite{car2} and in~\cite{car3}, she derived an asymptotic formula for the number of triples $(P_1,P_2,P_3)$ of irreducible polynomials in $\Fq[t]$ of equal degree which satisfy $P_1=P_2+P_3$.

Using sieve methods, Cherly~\cite{cherly} established that in $\Fq[t]$ with $q>2$ fixed, every polynomial $F$ of sufficiently high degree can be expressed as the sum of two polynomials of unequal degrees not larger than the degree of $F$ and each having at most four prime factors. 

Proving the analogon to Vinogradov's result, Effinger and Hayes~\cite{eff-hayes} showed that with the exception of polynomials of the form $t^{2}+a$ in characteristic~$2$, every noneven monic polynomial in $\Fq[t]$ of degree $n$ at least $2$ is the sum of three irreducible monic polynomials, one of degree $n$ and the other two of degree strictly smaller than $n$.

Hayes~\cite{hayes} gave a succinct proof of the fact that for $q$ sufficiently large relative to the degree $n$ of a polynomial $F\in \Fq[t]$, it is possible to write $F$ as the sum of two irreducibles both of degree $n+1$. Furthermore, he proved an asymptotic formula with $q\rightarrow\infty$ for the number of ways a given monic polynomial $F\in \Fq[t]$ of degree $n$ can be written as the sum of two irreducible monic polynomials both of degree $n+1$, assuming that $F$ is squarefree or that $\mathrm{char}(\Fq)\nmid n+1$. 

The most recent work on the Goldbach problem in the polynomial setting was done by Paul Pollack~\cite{pollack-thesis}~\cite{pollack-article}, who considered representations of a polynomial in $\Fq[t]$ of degree $n$ as the sum of two irreducibles of degree up to $n$. He sharpened the bound in~\cite{car1} for the number of polynomials which cannot be represented in this way. He also stated a conjecture on the number of such representations and proved it under the assumptions that $\gcd(q,2n)=1$ and that $q$ is large enough relative to $n$.

\bigskip
{\bf Acknowledgements}: I thank Olivier Wittenberg for comments clarifying several points in a preliminary draft of this article, the need to check that the rational maps $\beta_i$ are morphisms being one of them. I thank Gian Pietro Pirola for useful discussions and particularly for his simplification of my proof of lemma~\ref{le:irreducibility}. I also thank Najmuddin Fakhruddin and Donghoon Hyeon for useful discussions and Jean-Marie de Koninck and Claude Levesque for the invitation to the Qu\'ebec-Maine Conference on Number Theory and Related Topics $2006$, which got me started thinking about this problem. Special thanks are due to the anonymous referees for very good suggestions which lead to a corrected and much more readable final version. 
A large part of the research for this article was carried out while I was a member of the Korea Institute for Advanced Study. I am very grateful for the excellent working conditions which the institute provided to me. 
\vspace{1mm}

\noindent\emph{Notation.}
We write $\Falg$ for an algebraic closure of $\Fq$ and let $\F{q^{s}}$ denote its unique subfield of cardinality $q^{s}$. 
We say that a point on a plane curve is of order $n$ with some $n\geq 2$ if the scheme-theoretic intersection of the curve and its tangent at that point is of multiplicity $n$.

\section{Preparations for the proof of the main result}

Let $\mathfrak X$ denote the family of all algebraic curves over $\Falg$ which can be embedded into the affine plane $\A^{2}$ as curves of total degree $n-1$. 
We fix such an embedding of $\mathfrak X$ by a map
$$\begin{array}{rcl}
{\mathfrak X}      & \longrightarrow  & \A^{I}\times \A^{2},\\
{\mathfrak x}       & \mapsto               & (c)\times \{f_{1(c)}(x,t)=0\},
\end{array}$$
where $(c)$ is the coefficient vector of $f_{1(c)}$ belonging to an implicit choice of the ordering of monomials. 

We denote by $\mathfrak F$ the subset of $\mathfrak X$ for whose members the polynomial $f_{1(c)}(t+b,t)$ is monic and of degree $n-1$ in $t$ for any $b\in \Falg$. This means that for curves in $\mathfrak F$, the coefficients of the terms of total degree $n-1$ sum up to one, so the coefficient vector $(c)$ of such a curve is a point on an affine hyperplane $H\subset \A^{I}$. 

The dimension of the affine space parameterizing the polynomials of total degree $d$ in two variables is $(d+1)(d+2)/2$. For the degree $d=n-1$ of $f_1$, we use $H\cong \A^{I}$ with $I=n(n+1)/2-1$ as parameter space since we only consider normalized polynomials $f_1$ which correspond to points in the affine subspace $H$ of codimension 1. 

For any element $f_{1(c)}\in{\mathfrak F}$ with coordinates $(c)\in H_{\Falg}$, we set
$$
f_{2(c)}(x,t)  =  -f_{1(c)}(x,t)+F(t).
$$
We then consider the families of curves in $H_{\Falg}\times \P_{\Falg}^{2}$
\begin{equation}\label{equ:Ci}
\xymatrix@C=8ex{
C_{1(c)}\ar[dr]^{\alpha_1} & & C_{2(c)} \ar[dl]_{\alpha_2} \\
  &H_{\Falg},&    
}
\end{equation}
where $C_{i(c)}=\alpha^{-1}_i (c)$ denotes the Zariski closure in $\P^{2}_{\Falg}$ of the affine curve $f_{i(c)}(x,t)=0$. 

We shall define four Zariski open subsets $A_j\subset H$ which correspond to curves with various properties we shall need. 
We shall prove that every $A_j$ is a nonempty subset of $H$ and so $H\setminus A_j \neq H$ for $j=1,2,3,4$.
We use the fact that the embedding of $H\setminus A_j$ is contained in a hypersurface of degree bounded by, say, $d$. An upper bound on $\abs{H\setminus A_j}_{\Fq}$ is then provided by an upper bound on the number of $\Fq-$rational points on one hypersurface of degree $d$. Such a bound is established in the following lemma.

\begin{lemma}[IV \S 3 Lemma 3A in~\cite{schmidt}]\label{schmidt}
Let $f(x_1,\ldots, x_n)$ be a nonzero polynomial defined over $\Fq$ and of total degree $d$. Then $f$ has at most $dq^{n-1}$ zeros in $\Fq^{n}$. 
\end{lemma}

The hypersurfaces to which the previous lemma will be applied are often constructed using resultants, whose degrees are then of obvious importance. 

\begin{lemma}\label{le:resultant}
(a) The resultant $Res_x$ of two polynomials of degrees $d_1$ and $d_2$ in the variable $x$ with coefficients which are polynomials of total degree at most $B$ is a polynomial of total degree bounded by $B(d_1+d_2)$. The same is true for the projective resultant for which we homogenize the polynomials with respect to $x$. 

(b) Let $f_0, \ldots,f_k$ be homogeneous polynomials in the variables $x_0, \ldots, x_k$, all of total degree $d$. Then the resultant $R(f_0, \ldots,f_k)$ is a homogeneous polynomial of total degree $(k+1)d^k$ in the coefficients of each $f_i$. 
\end{lemma}

\begin{pf}
(a) The single-variable resultant is the determinant of a matrix of size $(d_1+d_2)\times (d_1+d_2)$ with entries which are polynomials of degree bounded by $B$. 

(b) This is a special case of proposition $2.3(ii)$ in~\cite{Jou} and also proposition $1.1$ of~\cite{GKZ} with the remark in the paragraph after its statement in chapter 13.   
\end{pf}

\bigskip
\begin{lemma}\label{le:smoothness}
Let the $C_{i(c)}$ be as defined in diagram~(\ref{equ:Ci}) and $A_1\subset H$ correspond to curves which are smooth. Then we have
$$\abs{H\setminus A_1}_{\Fq} \leq 3((n-1)^2+(n-2)^2) q^{I-1}.$$
\end{lemma}

\begin{pf}
Since a projective curve is complete, the projection to $H_{\Falg}$ of the singular loci of the two families $C_{i(c)}$ in $H_{\Falg}\times \P_{\Falg}^{2}$ is Zariski closed. 

For $p\nmid n-1$, we set $f_1=2x^{n-1}-t^{n-1}+d$ and then a short calculation shows that the corresponding $C_{1(c)}$ is smooth for any nonzero $d$ and that $C_{2(c)}$ is singular for at most $n-1$ different values of $d\in \Falg$. Likewise, for $p\mid n-1$ we take $f_1=2x^{n-1}+x-t^{n-1}+t^{n-2}$, for which both the resulting $C_{1(c)}$ and $C_{2(c)}$ are smooth. 
Therefore the points $(c)\in H$ whose fibres $C_{i(c)}$ are smooth curves form a nonempty Zariski open subset $A_1\subset H$.

We have to find an upper bound for the number of points $(c)\in H_{\Fq}$ such that there exists a singularity on one of the $C_{i(c)}$. The singularities of a projective curve $f(x,t,z)=0$ all satisfy the system of equations
$$
\begin{array}{l}
f_t:=\frac{df}{dt}=0,\\
f_x:=\frac{df}{dx}=0,\\
f_z:=\frac{df}{dz}=0.
\end{array}
$$
In the case of $p\mid n$ the equation $f(x,t,z)=0$ is a further restriction, but we ignore this slight improvement valid only if $p\mid n$. 

Let the degree of these three homogeneous polynomials be $d$. 
We now apply part $(b)$ of lemma~\ref{le:resultant} and deduce the existence of a polynomial $Res_{t,x,z}(f_t, f_x, f_z)$ in the coefficients of $\{f_t, f_x, f_z\}$ of total degree $3d^2$ with the property that $\{f_t, f_x, f_z\}$ have a common zero if and only if $Res_{t,x,z}(f_t, f_x, f_z)=0$. 

Applying this to the curves $C_{1(c)}$ and $C_{2(c)}$, lemma~\ref{schmidt} gives 
$$\abs{H\setminus A_1}_{\Fq} \leq 3((n-1)^2+(n-2)^2) q^{I-1}.$$
\end{pf}

\begin{lemma}\label{le:separability}
Let the $C_{i(c)}$ be as defined in diagram~(\ref{equ:Ci}) and $A_2\subset H$ correspond to curves with separable Gauss maps. Then we have
$$
\abs{H\setminus A_2}_{\Fq} \leq 16n^4 q^{I-1}.
$$
\end{lemma}

\begin{pf}
Under our assumption $p>2$, the Gauss map is separable if and only if the minimal intersection multiplicity of a tangent and the curve is equal to $2$~\cite[Prop. 4.4(b), with the notation $M(C)$ defined on page $1479$]{homma}. 

This condition of separability is satisfied for $C_{i(c)}$ if the curve has a simple tangent, that is if there exist $a,b,u,v\in \Falg$ and $h_{i(c)}(t)\in \Falg [t]$ such that we have
$$
\begin{array}{l}\label{tripletangent}
f_{i(c)}(at+b,t)=(u t-v)^{2}h_{i(c)}(t),\\
h_{i(c)}(v/u)\neq 0.
\end{array}
$$
The second equation defines an open subset of the closed set defined by the first equation in the variables $(c),a,b,u,v$ and the coefficients of $h_{i(c)}$. 
We project this open subset to the space $\A^{I}$ of coordinates of $f_i$ and restrict it to a subset of $H$. The intersection of the two subsets of $H$ arising from the two $C_{i(c)}$ we call $A_2$. 

In order to show that the projections of these open sets to $H$ are open, we need to argue that the family
$$
\{ f_{i(c)}(at+b,t)=(u t-v)^{2}h_{i(c)}(t)\} \rightarrow \A^{I}
$$
is flat over an open subset of $\A^{I}$. This is ensured by the generic flatness theorem of Grothendieck (\cite[Th\'{e}or\`{e}me 6.9.1]{AG-EGA4} or \cite[Theorem 14.4]{eisenbud-CA}). That a flat morphism is open is the contents of lemma $5.21(1)$ in~\cite{fantechi}. 

Therefore we have a Zariski open subset $A_2\subset H$ for whose points $(c)$ there exist simple tangents to the $C_{i(c)}$. 

In order to show that $A_2$ is nonempty, we let $F(t)=t^n+a_{n-1}t^{n-1}+\ldots +a_0$ and set 
$$f_{1}(x,t)=2x^{n-1}-x^{2}t^{n-3}+a_0+a_1 t+(a_2-1)t^2.$$
The case $n=3$ being obvious, for $n>3$ the curve $C_1$ and the line at infinity intersect in a point of order $2$, as do the curve $C_2$ and the line $x=t$. 

We need this argument only to deal with the case $p\mid n(n-1)(n-2)$. 
If $p\nmid n(n-1)(n-2)$, separability of the Gauss maps for both $C_{i(c)}$ with $(c)\in A_1$ is proved already in the last paragraph of the proof of proposition $3.1$ in~\cite{b-w}.

By the first sentence of the proof, curves with inseparable Gauss maps are excluded by the requirement to have at least one simple tangent. That all tangents to the curve $C_{2(c)}$ of the form $x=t+b$ are simple and that there is at least one such tangent is condition $(2)$ in lemma~\ref{le:spectangents} below. We can therefore use the same estimates as the ones derived in lemma~\ref{le:spectangents}:
$$
S_{2}(a_j) \leq n(2n-1)(2n(2n-1)-1)<8n^4 \mbox{ and }S_{3}(a_j) \leq n(2n-1). 
$$
For $C_{1(c)}$, the same argument involving condition $(1)$ and the analogous extension involving the existence of at least one tangent of the form $x=t+b$ gives
$$
S_{1} \leq n(2n-3)(2n(2n-3)-1)<8n^4 \mbox{ and }S_{4}(a_j) \leq n(2n-1). 
$$
Together this gives the sufficient bound
$$
\abs{H\setminus A_2}_{\Fq} \leq 16n^4 q^{I-1}.
$$
\end{pf}

\bigskip
We let $\beta_i$ denote the projection $C_{i(c)}\rightarrow \P^1$ from $M=(1,1,0)$ in homogenized coordinates $(x,t,z)$. 
We need to show that the rational maps $\beta_1$ and $\beta_2$ are in fact morphisms and this amounts to checking that the point $M$ does not lie on $C_{1(c)}$ or $C_{2(c)}$. We assume that $(c)\in H$ and monicity of the two polynomials $f_{i(c)}(t+b,t)$ for any $b$ implies both conditions $M\not\in C_{i(c)}, i=1,2$. 

\smallskip
If $k$ is a field, a finite $k$-scheme $X$ is said to have \emph{at most one double point} if
$n(X) \geq r(X)-1$, where $r(X)$ denotes the rank and $n(X)$ the geometric number of points of~$X$.

\smallskip
\begin{defn}(adapted from \cite{hurwitz})
Let $S_k$ be a Zariski open and dense subset of $\P^{1}_{k}$. A finite morphism $f\colon C \rightarrow S_{k}$ is called generic over $S_k$ if $f^{-1}(x)$ has at most one double point for all $x\in S_{k}$.
\end{defn}

\begin{lemma}\label{le:spectangents}
Let the $C_{i(c)}$ be as defined in diagram~(\ref{equ:Ci}). Assume the Gauss maps of both $C_{i(c)}$ are separable. Let $A_3$ denote the set of parameters whose corresponding curves satisfy the conditions

\begin{enumerate}
\item The morphism $\beta_1$ is generic over $\P^{1}$.
\item The morphism $\beta_2$ is generic over an affine subset $\A^{1}\subset \P^{1}$ and is ramified over at least one point of that $\A^{1}$. 
\item No line $x=t+b$ is tangent to both $C_{1(c)}$ and $C_{2(c)}$.
\item The line at infinity is not tangent to $C_{1(c)}$. 
\end{enumerate}

Then we have
$$\abs{H\setminus A_3}_{\Fq} \leq 40 n^4 q^{I-1}.$$
\end{lemma}

\begin{pf}
We denote the dual projective plane by $\P^{2\star}$. Let the algebraic set $X$ consist of the tangents which intersect $C_{1(c)}\cup C_{2(c)}$ with $(c)\in A_1\cap A_2$ in more than one double point and which we shall call special tangents. The possible special tangents have intersection multiplicity larger than $2$ with one of the curves $C_{i(c)}$, they are bitangent to one of the $C_{i(c)}$ or are tangent to both $C_{i(c)}$. 

We projectivize the tangents' coordinates $a,b$ in $\A^{2\star}$ to $\P^{2\star}$ and then we have the projections
$$
\xymatrix@C=8ex{
X\subset A_1\cap A_2\times \P^{2\star}\ar[d]^{\gamma}\ar[r]^{\hspace{12mm}\delta} &\P^{2\star}\\
A_1\cap A_2
}
$$

We now consider the Zariski closed set $X_1\subset \gamma (X)$ with $(c)\in X_1$ if there exists a special tangent whose affine equation is $x=t+b$. Note that for any $(c)$, the line at infinity intersects the curve $C_{2(c)}$ in a point of multiplicity $n$. 

We consider the curves $C_{1(c)}$ given by an equation of the form
$$
f_{1(c)}=t^{n-1}+g_{1(c)}(x,t),
$$
where $g_{1(c)}(x,t)$ is of degree at most $n-2$. We let $H_3\subset H$ denote the subspace of $H$ corresponding to curves of this form. 

\smallskip
In the first paragraph of the proof of proposition 3.1 in~\cite{b-w}, it is pointed out that the proof of the proposition shows a stronger result, namely that for all $(c)\in A_1\cap A_2$, there are only finitely many points in $\delta(\{(c)\times \P^{2\star}\}\cap X)$. 
Therefore a separable Gauss map implies that the curve has only finitely many special tangents, so we can choose an $r\in \Falg$ such that in the new coordinates $(x',t)$ with $x'=rx$, the coefficients $(c)\in H_3$ are transformed into $(c')\in H_3$ for whose associated curves $C_{i(c')}$ there is no special tangent with affine equation $x=t+b$. It follows that $X_1\neq \gamma (X)$ and the form of the equations of curves in $H_3$ shows that the polynomials defining the curves $C_{1(c')}$ are monic. 

The line at infinity is not tangent to $C_{1(c)}$ for $(c)$ in an open subset of $H$ and the choice $f_1=2x^{n-1-j}t^{j}-t^{n-1}$ with $j=1$ for $p\mid n-1$ and $j=0$ otherwise shows that this open subset is nonempty.  

We have ramification of $\beta_2$ over $\A^1$ if the resultant 
$$
Res_t (f_2(t+b,t), \frac{d}{dt}f_2(t+b,t))
$$
is a polynomial of nonzero degree in the variable $b$. This condition defines a Zariski open subset of $\A^{I}\times \A^1$, where the variable of $\A^1$ is $b$.  We project this set to $\A^I$ and get a Zariski open subset of curves $C_{2(c)}$ which ramify over $\A^1$. 

We are left to show that this Zariski open subset is actually nonempty. As before we let $F(t)=t^n+a_{n-1}t^{n-1}+\ldots +a_0$ and we choose $f_1(x,t)=t^{n-1}+a_{n-2}t^{n-2}+\ldots +a_0+x-t-c_2 t^2$ for $n\geq 3$ and $f_1(x,t)=2x-t$ for $n=2$. In both cases it is straightforward to check that a suitable choice of $b$ and $c_2$ will ensure that $f_2(t+b, t)$ has a root of multiplicity two. 

We have shown that there exists a nonempty open subset $A_3\subset H$ such that all $f_{i(c)}$ with $(c)\in A_3$ satisfy the following conditions: the associated morphism $\beta_1$ is generic over $\P^1$, the associated morphism $\beta_2$ is generic over an affine subspace $\A^1\subset \P^1$ and ramifies over $\A^1$, no line $x=t+b$ is tangent to both $C_{i(c)}$ and the line at infinity is not tangent to $C_{1(c)}$.

\smallskip
We proceed to estimate the sizes of the Zariski closed subsets to be excluded from $H$. 

Condition $1$ of genericity of the morphism $\beta_1$: 

We projectivize the affine curve $f_{1(c)}(t+b,t)=0$ only with respect to the variable $t$. This produces a polynomial $F_{1(c)}(t+b,t,z)=0$ homogeneous in the variables $t$ and $z$.

The line at infinity is treated under condition $4$. We are left with special tangents given by equations $x=t+b$. Values $b$ such that $x=t+b$ is tangent to the projective closure of the affine curve $f_{1(c)}(x,t)=0$ are solutions of 
$$R_{1}(a_j, b)=Res_{t,z}(F_{1(c)}(t+b,t,z), \frac{d}{dt} F_{1(c)}(t+b,t,z))=0,$$ 
where $Res_{t,z}$ stands for the projective resultant with respect to the variables $t$ and $z$. 
There are only simple tangents of the form $x=t+b$ to the projective closure of $f_{1(c)}(x,t)=0$ if $R_{1}(a_j, b)$ does not have repeated roots as a polynomial in $b$, a condition equivalent to 
\begin{equation}\label{eq:multtangents1}
S_{1}(a_j) = Res_b(R_{1}(a_j, b), \frac{d}{db} R_{1}(a_j, b))\neq 0.
\end{equation}
Here also the $a_i$ enter the calculation of the total degree. Hence, using part $(a)$ of lemma~\ref{le:resultant}, for the estimates of the total degrees of $R_{1}$ and $S_{1}$ we get 
$$
\deg R_{1}\leq n(2n-3)\mbox{ and }\deg S_{1} \leq n(2n-3)(2n(2n-3)-1).
$$ 

\smallskip
Condition $2$ of genericity and ramification of the morphism $\beta_2$:

Here we use the same approach as for $C_{1(c)}$, the only difference being that the line at infinity is always tangent to $f_{2(c)}$. 
The polynomial $f_{2(c)}$ is monic in $t$ and so the bounds work out to 
$$
\deg R_{2}\leq n(2n-1),
$$
$$
\deg S_{2} \leq n(2n-1)(2n(2n-1)-1).$$

The morphism $\beta_2$ is ramified over $\A^1$ if 
$$R_3(a_j, b)=Res_t (f_2, \frac{d}{dt} f_2)$$
regarded as a polynomial in $b$ is not a constant. For this it suffices to demand that the coefficient of one term $b^k$ with $k\neq 0$ be nonzero. Since $R_3(a_j, b)$ is of total degree $n(2n-1)$, we use this as a bound on the degree of the coefficient of $b^k$ and get  
$$\deg S_3(a_j)\leq n(2n-1).$$ 

\smallskip
Condition $3$ on lines tangent to both $C_{i(c)}$: 

The line $x=t+b$ is tangent to the projective closure $C_{i(c)}$ if 
$$R_i(a_j, b)=Res_{t,z}(F_i(t+b,t,z), \frac{d}{dt} F_i(t+b,t,z))=0$$ 
and so the line is tangent to both $C_{i(c)}$ if $Res_b(R_1, R_2)=0$. An estimate derived in analogy to condition $1$ works out to 
$$\deg Res_b(R_1, R_2) \leq n(2n-1)n(4n-4)=4n^2(n-1)(2n-1).$$
The line at infinity will again be treated as part of condition $4$. 

\smallskip
Condition $4$ on the line at infinity: 

This condition is equivalent to the polynomial consisting of the monomials of degree $n-1$ of $f_1$ not having a double root. This corresponds to a resultant of size $2n-3$ with coefficients of degree $1$, so we have to remove one hypersurface of degree $2n-3$ from $H$.

In total we have five degrees which are all bounded by $8n^4$, so a sufficient bound is
$$\abs{H\setminus A_3}_{\Fq} \leq 40 n^4 q^{I-1}.$$  

\end{pf}

\bigskip
As for notation, we set $k_1=k(a_i, b)$ and let $r$ be the image of $t$ in the quotient ring $K:=k_1[t]/f_{2(c)}(t)$.

\begin{lemma}\label{le:irreducibility}
Let the $C_{i(c)}$ be as defined in diagram~(\ref{equ:Ci}). Let $A_4$ denote the set of parameters whose corresponding curves satisfy the condition that in the factorization of $f_{2(c)}$ as 
$$
f_{2(c)}(t)=(t-r)g(t),
$$
the polynomial $g(t)$ is irreducible over the field $k_1[r]$. 

Then we have
$$\abs{H\setminus A_4}_{\Fq} \leq (16^{n^{4}}n^4 (n+1)^{8n^4+1}) q^{I-1}.$$

\end{lemma}

\begin{pf}
We consider the family of global fields generated by roots of $f_{2(c)}(t)=0$ and parameterized by $(c)$. 
Consider $f_{2(c)}(t+b,t)$ a polynomial in the variable $t$ with coefficients in $k_1$ and factorize it as in 
$$
f_{2(c)}(t)=(t-r)g(t).
$$
Since $f_{2(c)}$ is monic, the element $r$ and all its conjugates are integral in the extension field $K$, so in a factorization
$$
f_{2(c)}(t)=(t-r)g_1(t)g_2(t),
$$
all coefficients of the $g_i(t)$ are also integral. 

Choosing any values of the $a_i$, the set $\{1,r,\ldots ,r^{n-1}\}$ is a basis for $K$ over $k(b)$ consisting of integral elements. The discriminant of $(1,r,\ldots ,r^{n-1})$ we denote by $d$. We recall the standard fact from algebraic number theory that every integral element in the field $K$ is a linear combination of $\{ \frac{1}{d}, \ldots , \frac{r^{n-1}}{d}\}$ with coefficients in $k[b]$. 
Assume $f_{2(c)}(t)$ factorizes as 
\begin{equation}\label{equ:factorize}
f_{2(c)}=(t-r)(t^m + u_{m-1}r^{m-1}+\ldots )(t^{n-m-1}+v_{n-m-2}t^{n-m-2}+\ldots ).
\end{equation}  
By the previous paragraph, the $u_i, v_i$ are integral and can be chosen to have degree at most $n-1$ in $r$. 
All these integral coefficients $u_i, v_i$ are then of the form 
$$
\frac{1}{d}\sum_{j=0}^{n-1} s_j r^{j}.
$$
We want to derive a bound for the degree of the polynomials $s_j\in k[b]$. The absolute value $\abs{s(b)}_{\infty}:=q^{\deg_b(s)}$ on $k[b]$ can be extended to a finite extension of $k(b)$. The extension which is relevant to our purposes is the splitting field of $f_{2(c)}(t)$ with roots $\alpha_i$. 
The product of all roots  $\alpha_i$ is equal to the constant term of $f_{2(c)}(t)$ as a polynomial in $t$, so we have 
$$
\abs{\prod \alpha_{i}}_{\infty}\leq q^n
$$
and therefore all coefficients $u_i, v_i$ have absolute value bounded by $q^n$. 
We therefore get
\begin{equation}
\label{equ:bounddiscr}
\deg_b(\frac{s_j}{d})\leq n.
\end{equation}

Now we need an estimate for the degree of $d$. 
As an expression for $d$ in terms of the resultant, we have $d=Res(f_2(t), f_2'(t))$\cite[Prop. IV.2.7]{Lorenzini} and this is the determinant of a matrix of size $2n-1$. By part $(a)$ of lemma~\ref{le:resultant}, the total degree is therefore $\deg(d)\leq (n+1)(2n-1)$ and $\deg_b(d)\leq n(2n-1)$. 
The bound~(\ref{equ:bounddiscr}) then implies $\deg_b (s_j)\leq n+n(2n-1)=2n^2$. 
One of the coefficients $u_i, v_i$ needs at most $2n^3$ coefficients and one complete factorization, with $n-1$ such coefficients, needs $v:=2n^3(n-1)$.
Comparing coefficients in $k$, we have $v+1$ homogeneous equations $e_i=0$ in $v+1$ variables. 
 
Degree estimate of the equations:
Multiplying out equation~(\ref{equ:factorize}), the right hand side has coefficients in the form of terms
$$
\frac{u_i}{d}+\frac{u_i v_j}{d^2}.
$$
In the numerator, we can apply the quotient operation of the ring $K$ in the form of subtraction of suitable multiples of $f(r)$ by powers of $r$ and polynomials in $b$ to get representatives of degree at most $n-1$ in $r$. 
Multiplying by $d^2$, the highest degree terms are on the left hand side: the total degree is bounded by
$$
n+1+(n+1)^2(2n-1)^2\leq 4(n+1)^4.
$$

The polynomial $g(t)$ is reducible only if it is divisible by a factor of degree up to at most $(n+1)/2$, which corresponds to at most $(n+1)/2$ equations of the form 
\begin{equation}\label{irred}
f_{2(c)}(t)=(t-r)g_1(t)g_2(t). 
\end{equation}

Eliminating the new $v+1$ variables in the equations~(\ref{irred}), it follows that the set $A_4$ of $(c)$ with $f_{1(c)}$ inducing $f_{2(c)}$ which satisfy our condition is Zariski open. Choosing an $f_1(t,b)$ such that $f_2(t,b)=t g(t)$ with $g(t)=t^{n-1}+g_1(t,b)$ irreducible and with $\deg g_1 \leq n-2$ shows that $A_4$ is nonempty. More generally, we can argue that for $f_{2(c)}(t)$ a polynomial with Galois group the symmetric group $S_n$, the resulting $v+1$ equations in $v+1$ homogeneous variables have no solution and so 
$$
Res(e_0, e_1,\ldots ,e_{v})
$$
is a nonvanishing polynomial in the coefficients $a_i$.  

We apply lemma~\ref{le:resultant} to $v+1$ equations of total degree at most $4(n+1)^4$. The degree of the resultant corresponding to one factorization of $f_{2(c)}$ is therefore bounded by $2 n^{4}16^{n^4} (n+1)^{8n^4}$. 

For $(n+1)/2$ factorizations, we therefore have
$$\abs{H\setminus A_4}_{\Fq} \leq (16^{n^4} n^{4} (n+1)^{8n^4+1}) q^{I-1}.$$

\end{pf}

\section{Proof of Theorem~\ref{cor}}
We follow the proof of theorem $1.1$ in~\cite{b-w}, adding necessary adaptations and explicit estimates for lower bounds of $q$. The notations and conventions of the previous sections are used freely in the proof. 

For a rough explanation of the main ideas, we quote the original version of the theorem whose proof we are going to adapt. 

\begin{thm}[Theorem 1.1 in~\cite{b-w}]
\label{maintool}
Let $\Fq$ be a finite field of characteristic $p$ and cardinality $q$.
Let $\tuple{f_1}{f_n} \in \Fq[t,x]$ be irreducible polynomials whose
total degrees $\deg(f_i)$ satisfy $p \nmid \deg(f_i) (\deg(f_i)-1)$
for all~$i$. Assume that the curves $C_i \subseteq \P^2_\Fq$ defined
as the Zariski closures of the affine curves
$$
f_i(x,t)=0
$$
are smooth. Then, for any sufficiently
large $s \in \N$, there exist \mbox{$a, b \in \F{q^s}$} such that
the polynomials $\tuple{f_1(at+b,t)}{f_n(at+b,t)} \in \F{q^s}[t]$
are all irreducible.
\end{thm}

\smallskip
Let $F(t)$ be a polynomial in $\Fq[t]$ and suppose theorem~\ref{maintool} can be applied to 
$$\begin{array}{l}
f_1(x,t),\\
f_2(x,t)=-f_1(x,t)+F(t)
\end{array}$$
with some $f_1$ to be chosen. In view of $F=f_1+(-f_1+F)$, this would be enough to represent $F(t)$ as the sum of two irreducibles if $q$ is large enough. 

Key to proving the existence of such an $f_1$ is the observation that theorem~\ref{maintool} demands genericity conditions on the polynomials $f_i$. The condition involving indivisibility by $p$ is assumed to ensure separability of the Gauss maps of the $C_i$, which is a genericity condition as well.
For each one of these conditions we shall show that there exists a nonempty Zariski open subset of $H_{\Falg}$ whose points $(c)$ satisfy it. 

In the proof of theorem~\ref{maintool}, the intersection of the line $x=at+b$ and a curve $C_i$ is  viewed as the fibre of a projection $C_i\rightarrow \P^1$ from a point $M$ in $\P^2$. For the purpose of that proof, it turns out that $M$ can be chosen freely in a nonempty Zariski open subset of $\P^2$. In order to ensure that the resulting irreducible polynomials are monic, we need to fix the value of $a$ in advance. The condition $a=1$ is equivalent to choosing $M$ as the point at infinity with $x=1$ and $t=1$. This choice of $M$ imposes some further necessary conditions on the two polynomials $f_i$. 

More precisely, from the number of points $(c)$ in $H_{\Fq}$ parameterizing the polynomials $f_{1(c)}$ and $f_{2(c)}$, we shall deduct an upper bound for the number of $\Fq-$rational points in $H\setminus A_j$ for $j=1,2,3,4$. 
The lower bound on $q$ stated in the theorem will then imply that the remaining set of points $(c)$, and therefore the set of usable polynomials $f_{1(c)}$ and $f_{2(c)}$, is nonempty. Furthermore, we shall show that the lower bound on $q$ is also large enough for the \v{C}ebotarev Density Theorem in the form of theorem $3.5$ in~\cite{b-w} to be applicable to the field extensions defined by the resulting $f_{1(c)}$ and $f_{2(c)}$. 

Subject only to a lower bound on $q$, the adapted argument then produces an element $b_0\in \F{q}$ such that both $f_1(t+b_0,t)$ and $f_2(t+b_0,t)$ are monic and irreducible in $\F{q}[t]$, which is enough to prove theorem~\ref{cor}. 

\bigskip
We let~$\Falg$ denote an algebraic closure of~$\Fq$ and~$M \in \P^2(\Fq)$
denote the point at infinity with coordinates $x=1$, $t=1$. Furthermore, we use the letter $k$ to denote the field $\Fq$. 

We have proved four lemmata to the effect that if $q$ is large enough, there exist curves $C_1$ and $C_2$ with certain properties. Lemma~\ref{le:smoothness} ensures smoothness, lemma~\ref{le:separability} separability of the Gauss map, lemma~\ref{le:spectangents} excludes existence of undesired special tangents and lemma~\ref{le:irreducibility} deals with irreducibility of $f_2$ in an extension field.

We use our assumption on the size of $q$ and note that it is three times the largest lower bound demanded by the four lemmata listed in the previous paragraph, where we have combined the bounds in lemma~\ref{le:separability} and lemma~\ref{le:spectangents}. We can therefore assume that all conclusions of the four lemmata hold for the chosen $f_1$ and $f_2$. 

The following setup for the proof is almost the same as the one used for the main result in~\cite{b-w}. 
We denote by $\beta_i \colon (C_i)_k \rightarrow \P^1_k$ the
$k$-morphism $\P^2_k \setminus \{M\} \rightarrow \P^1_k$
defined by projection from~$M$. Furthermore, we write $\kappa(X)$ for the
function field of $X$ when $X$ is an integral scheme.
The morphism $\beta_1$ is finite of degree $n-1$ and is
generic over $\A^1$; this is the contents of part 1 of lemma~\ref{le:spectangents}. The morphism $\beta_2$ is finite of degree $n$ and is
generic over $\P^1$ and this is proved in part 2 of lemma~\ref{le:spectangents}. 
 Being generic over a nonempty open subset of $\P^1$, both morphisms $\beta_i$ are separable;
therefore there exists a smooth, complete,
connected curve $C'_i$ over~$k$ and a finite morphism $C'_i \rightarrow (C_i)_k$, such that
the induced field extension $\kappa(C'_i)/\kappa(\P^1_k)$ is a Galois closure of
$\kappa((C_i)_k)/\kappa(\P^1_k)$.
We introduce the notation $K=\kappa(\P^1_k)$,
$K_i=\kappa(C'_i)$, $G_i=\Gal(K_i/K)$ and $H_i=\Gal(K_i/\kappa((C_i)_k))$.

\smallskip
The following proposition is enough to deal with the case of $\beta_{2}$. 
\begin{prop}
\label{prelproporig}
Let $C$ be a regular, complete, geometrically irreducible curve over a field~$k$,
endowed with a finite separable generic morphism $f \colon C \rightarrow \P^1_k$.
Let $C'$ be a regular, complete, irreducible curve over $k$, and $g \colon C' \rightarrow C$
be a finite morphism. Assume that the finite extension $\kappa(C')/\kappa(\P^1_k)$ is a Galois
closure of the subextension $\kappa(C)/\kappa(\P^1_k)$.
We denote respectively by $G$ and $H$ the Galois groups of
$\kappa(C')/\kappa(\P^1_k)$ and $\kappa(C')/\kappa(C)$.
Then $C'$ is geometrically irreducible over~$k$ and the morphism
$$
G \longrightarrow \symm{H \backslash G}
$$
induced by right multiplication is an isomorphism.
Moreover, all the ramification indices of $\kappa(C')/\kappa(\P^1_k)$
are at most $2$.
\end{prop}

\begin{pf}
Proposition $2.2$ in~\cite{b-w}. 
\end{pf}

\smallskip
For the polynomial $f_{1(c)}$, we need the following adapted version of the above proposition. In order to formulate it, we introduce the following notation. 

Let $k'$ denote the algebraic closure of $k$ in $\kappa(C')$.
We denote respectively by $G'$ and $H'$ the subgroups of $G$ defined
by the subfields $\kappa(\P^1_{k'})$ and $\kappa(C_{k'})$ of $\kappa(C')$, so that
we have a canonical commutative diagram as follows, where the labels indicate the Galois
groups of the generic fibres:
$$
\xymatrix@C=8ex{
C_{k'} \ar[dd] \ar[rr] && \P^1_{k'} \ar[dd] \\
& C' \ar[ul]_{H'} \ar[ur]^{G'} \ar[dl]_H \ar[dr]^G \\
C \ar[rr] && \P^1_k
}
$$
Let us endow $H \backslash G$ (resp.~$H' \backslash G'$) with the action of~$G$ (resp.~$G'$)
by right multiplication.

\begin{prop}
\label{prelpropnew}
Let $C$ be a regular, complete, geometrically irreducible curve over a field~$k$ with an embedding of degree $n$ in $\P^2$,
endowed with a finite separable morphism $f \colon C \rightarrow \P^1_k$ which is generic over $\P^1_k\setminus \{ P\}$, ramified over $\P^1\setminus \{P\}$ and whose fibre above $P$ is a point of multiplicity $n$. We assume that in the factorization
\begin{equation}\label{equ:irred}
f(t)=(t-r)g(t),
\end{equation}
the factor $g(t)$ is irreducible over $k(x,r)$. 
Let $C'$ be a regular, complete, irreducible curve over $k$, and $g \colon C' \rightarrow C$
be a finite morphism. Assume that the finite extension $\kappa(C')/\kappa(\P^1_k)$ is a Galois
closure of the subextension $\kappa(C)/\kappa(\P^1_k)$.
We denote respectively by $G$ and $H$ the Galois groups of
$\kappa(C')/\kappa(\P^1_k)$ and $\kappa(C')/\kappa(C)$.
Then $C'$ is geometrically irreducible over~$k$ and the morphism
$$
G \longrightarrow \symm{H \backslash G}
$$
induced by right multiplication is an isomorphism.
Moreover, all the ramification indices of $\kappa(C')/\kappa(\P^1_k)$ over $\P^1_k\setminus \{ P\}$
are at most $2$.
\end{prop}

\begin{pf}
Every nontrivial inertia group of $G'$ over $\P^1_k\setminus \{ P\}$ has order $2$ and acts by a transposition on $H' \backslash G'$: this is proved as in proposition~\ref{prelproporig}. The assumption on the factor $g(t)$ in equation~(\ref{equ:irred}) implies that $G'$ is doubly transitive acting on $H' \backslash G'$. It follows immediately from lemma $4.4.3$ in~\cite{serre} that a doubly transitive group containing a transposition is the full symmetric group, hence the proposition.

\end{pf}

\bigskip
The propositions~\ref{prelproporig} and~\ref{prelpropnew} together with lemma~\ref{le:irreducibility} now show that for both indices $i=1,2$, the curve $C'_i$ is
geometrically connected over~$k$ and
the group $G_i$ is canonically isomorphic to $\symm{H_i \backslash G_i}$. 

Let $R_i \subset \P^1_k$ denote the branch locus of the morphism $C'_i \rightarrow \P^1_k$.

\begin{prop}
\label{prop:ridisjoint}
The subsets $R_i \subset \P^1_k$ for $i=1,2$ are disjoint.
\end{prop}

\begin{pf}
We shall need the following well-known lemma, which is a direct consequence
of Lemma $5.1$ in~\cite{b-w}.

\begin{lemma}
Let $E/K$ be a finite separable extension of global fields, and let $L$ be a Galois closure of $E/K$. Then a finite place of $K$ is unramified in $E$ if and only if it is unramified in $L$. 
\end{lemma}

The lemma shows that $R_i$ is also the branch locus of the morphism
$(C_i)_k \rightarrow \P^1_k$. An $\Falg$-point of $R_i \cap R_j$ therefore gives rise
to a line in $\P^2_\Falg$ which is both tangent to $(C_i)_\Falg$ and $(C_j)_\Falg$
and contains $M$, which is ruled out by lemma~\ref{le:spectangents}. 
\end{pf}

\bigskip
Let $L$ denote the ring $K_1 \otimes_K K_2$.

\begin{prop}
The ring $L$ is a field, and $k$ is separably closed in~$L$.
\end{prop}
\begin{pf}
This is proposition $3.4$ in~\cite{b-w} for $n=2$. 
\end{pf}

\bigskip
The following two paragraphs are quoted with small adaptations from~\cite{b-w}. 
Let $C'$ denote a smooth complete connected curve over~$k$ with function field $L$.
There is a natural finite morphism $\psi \colon C' \rightarrow \P^1_k$, which is
generically Galois and therefore separable.
We denote by $g$ the genus of $C'$, by $G$ the group $\Gal(L/K)$, by~$N$ the
degree of $\psi$, and by $(x, L/K)$ the Artin symbol of the extension $L/K$ above
a closed point $x \in \P^1_k$ which does not ramify in $L$.
We would now like to find a rational point of $\P^1_k$ above which the fibre of $\psi$ is
integral. To this end, we resort to an effective version of the \v{C}ebotarev
theorem for function fields, due to Geyer and Jarden. The following is a weak consequence
of~\cite[Proposition~13.4]{gejapaper}.

\begin{thm}
\label{gejathm}
Let $c$ be a conjugacy class in $G$. We denote by $P(L/K,c)$ the set of rational
points $x \in \P^1(k)$ outside the branch locus of $C' \rightarrow \P^1_k$
such that \mbox{$c=(x,L/K)$}.
Then one has
\begin{equation}
\label{gejaeq}
\abs{P(L/K,c)} \geq \frac{1}{N}\left( q^s - (N+2g)q^{s/2} - N q^{s/4} - 2(g+N) \right) \text{.}
\end{equation}
\end{thm}

Some preparation is in order before applying Theorem~\ref{gejathm}: to be able to deduce
from it that $P(L/K,c)$ is non-empty as soon as $s$ is chosen large enough, we need to make
sure that the right-hand side of~(\ref{gejaeq}) does grow when $s$ goes to infinity. For instance,
it suffices to establish that $N$ and $g$ are bounded independently of $M$ and $s$.
The integer~$N$ is obviously independent of the choices made: it is equal
to $\deg(f_{1(c)})! \deg(f_{2(c)})!=(n-1)! n!$. We shall actually prove that $g$ is also independent
of~$s$.

Remark: In~\cite{b-w} this is achieved by an application of the Riemann-Hurwitz formula, but in the case of wild ramification, i.e., if $p$ divides the degree of one of the polynomials, the genus is more difficult to bound in terms of ramification data alone. We therefore choose to rely on the total degrees for a bound on $g$. 

\begin{prop}\label{prop:genusbound}
Consider two curves
$$
\begin{array}{l}
f_1(x,t)=0\\
f_2(x,t)=0
\end{array}
$$
with $f_1$ of total degree $n-1$ and $f_2$ of degree $n$ and the Galois groups of both polynomials the full symmetric groups of size $n-1$ and $n$, respectively. We use the notation of function fields as above. Assume that $q> (n-1)^{(n-1)^2}n^{n^2}$. Then for the genus of $K_1\otimes K_2$ we have the upper bound
$$
g\leq (n!)^8. 
$$
\end{prop}

\begin{pf}
Proposition $3.4$ in~\cite{b-w} implies that the fields $K_1$ and $K_2$ are linearly disjoint over $K$. 
We consider the polynomials $f_i$ in the variable $x$ with coefficients in $k(t)$. We can construct $K_1\otimes K_2$ by adjoining all roots of the two $f_i$. 
By assumption, we have $Gal(K_1\otimes K_2)=S_{n-1}\times S_n$. We construct a simple extension $K_1\otimes K_2=K(\phi)$ by adjoining roots. 
Let $\alpha, \beta$ be roots of $f_1 f_2=0$. Every subfield of $K_1\otimes K_2$ corresponds to a subgroup of $S_{n-1}\times S_n$. The simplest bound for the number of such subgroups is $(n-1)^{(n-1)^2}n^{n^2}$, 
if $q>(n-1)^{(n-1)^2}n^{n^2}$, there exist $c,c'\in \Fq$ with
$$
K(\alpha+c\beta)=K(\alpha+c'\beta).
$$
It is an exercise concerning primitive elements in field theory that we then have $K(\alpha, \beta)=K(\alpha + c\beta)$. 

Let $v$ denote the degree valuation on $k[t]$ and its extension to the ring of integers of $K_1\otimes K_2$. We now embed $K(\phi)$ in the projective plane as
$$
\prod_{\gamma\in Gal}(x-\gamma(\phi))
$$
and it follows that the valuation $v$ of any coefficient is bounded by $n! (n-1)!n\leq (n!)^2$. The degree in $x$ is bounded by $n! (n-1)!$ and the degree in $t$ by $(n!)^2$ and so the total degree by $(n!)^4$. 
Singularities only decrease the genus, so we have
$$
g\leq \frac{(deg -1)(deg -2)}{2}\leq (n!)^8.
$$
\end{pf}

\bigskip
The assumption on the size of $q$ suffices for the application of proposition~\ref{prop:genusbound}. Then it is  almost immediate that the right-hand side of inequality~(\ref{gejaeq}) is $\geq 2$. 

The remaining part of the proof follows very closely the version in~\cite{b-w} and only needs insignificant adaptations. 
The canonical isomorphism $G = G_1 \times G_2 =
\symm{H_1 \backslash G_1} \times \symm{H_2 \backslash G_2}$ allows us to choose
an element $\sigma \in G$ whose projection 
in $G_i$ acts transitively on $H_i \backslash G_i$ for $i=1,2$.
Theorem~\ref{gejathm} now ensures the existence of a rational point $x \in \P^1(k)$
outside $R_1 \bigcup R_2$ and
such that $\sigma=(x,L/K)$. As the image of $(x,L/K)$ in $G_i$ is $(x,K_i/K)$,
it follows from Lemma 5.1 in~\cite{b-w} and the definition of $\sigma$ that
$\phi_i^{-1}(x)$ is irreducible. Moreover, the $k$-scheme $\phi_i^{-1}(x)$ is \'{e}tale since
$x \not\in R_i$, and hence it is integral.
Corresponding to the point $x$, there exists a $b_0 \in k$
such that for every $i=1,2$, the scheme
$\Spec(k[t]/(f_i(t+b_0,t)))$ is an open subscheme of $\phi_i^{-1}(x)$;
as the latter scheme is integral, the polynomials $f_1(t+b_0,t)$ and $f_2(t+b_0,t)$ must
be irreducible.  

Both $f_1(t+b_0,t)$ and $f_2(t+b_0,t)$ are monic and irreducible elements of $\Fq[t]$ with $\deg(f_1)=n-1$ and $\deg(f_2)=n$ and the proof is complete.

\bigskip
\noindent Andreas O. Bender\\
Pohang Mathematics Institute\\
POSTECH\\
San 31 Hyoja-dong\\
Pohang 790-784\\
Republic of Korea\\
\smallskip 

\noindent e-mail: \texttt{andreasobender@mac.com}\\
\texttt{http://www.andreasobender.org}

\begin{thebibliography}{11}

\bibitem[1]{CRnote}
\textsc{Andreas O. Bender}, \rm{Decompositions into sums of two irreducibles in $\Fq[t]$, C.R. Math. Acad. Sci. Paris \textbf{346} no. 17-18 (2008), 931--934, also available from \texttt{http://dx.doi.org/10.1016/j.crma.2008.07.025}}

\bibitem[2]{b-w}
\textsc{Andreas O. Bender, Olivier Wittenberg}, \rm{A potential analogue of Schinzel's hypothesis for polynomials with coefficients in $\Fq[t]$, Int. Math. Res. Not. \textbf{36} no. 36 (2005), 2237--2248, also available from \texttt{http://arxiv.org/abs/math/0412303}}

\bibitem[3]{car1}
\textsc{Mireille Car}, \rm{Le probl\`eme de Goldbach pour l'anneau des polyn\^omes sur un corps fini, C.R. Acad. Sci. Paris, S\'er. A-B \textbf{273} (1971), A201--A204.}

\bibitem[4]{car2}
\textsc{Mireille Car}, \rm{Le th\'eor\`eme de Chen pour $\Fq[X]$, Dissertationes Math. (Rozprawy Matematyczne) \textbf{223} (1984), Polska Akademia Nauk. Instytut Matematyczny.}

\bibitem[5]{car3}
\textsc{Mireille Car}, \rm{The generalized polynomial Goldbach problem, J. Number Theory \textbf{57} no. 1 (1996), 22--49.}

\bibitem[6]{chen}
\textsc{Jing-Run Chen}, \rm{On the representation of a larger even integer as the sum of a prime and the product of at most two primes, Sci. Sinica \textbf{16} no. 2 (1973), 157--176, also available in }\it{Yuan Wang,} \rm{ed., The Goldbach Conjecture, second edition, World Scientific Publishing, Singapore, 2002.}

\bibitem[7]{cherly}
\textsc{J\o rgen Cherly}, \rm{A lower bound theorem in $\Fq[x]$, J. Reine Angew. Math. \textbf{303/304} (1978), 253--264.}

\bibitem[8]{davenport}
\textsc{Harold Davenport}, \rm{Multiplicative number theory, Graduate Texts in Mathematics \textbf{74}, third edition, revised by Hugh L. Montgomery, Springer--Verlag, New York, NY, 2000.}

\bibitem[9]{eff-hayes}
\textsc{Gove W. Effinger, David R. Hayes}, \rm{A complete solution to the polynomial $3$-primes problem, Bull. Amer. Math. Soc. (N.S.) \textbf{24} no. 2 (1991), 363--369.}

\bibitem[10]{eff-hayes-book}
\textsc{Gove W. Effinger, David R. Hayes}, \rm{Additive number theory of polynomials over a finite field, Oxford University Press, New York, NY, 1991.}

\bibitem[11]{eisenbud-CA}
\textsc{David Eisenbud}, \rm{Commutative algebra with a view toward algebraic geometry, Graduate Texts in Mathematics \textbf{150}, Springer--Verlag, New York, NY, 1995.}

\bibitem[12]{euler}
\textsc{Leonhard Euler}, \rm{Letter to Christian Goldbach dated 30th June 1742,
Adolf P. Juskevic, Eduard Winter, eds., Leonhard Euler und Christian Goldbach, Briefwechsel 1729--1764, Akademie--Verlag, Berlin, 1965, also available as Lettre XLIV from \texttt{http://www.math.dartmouth.edu/$\sim$euler/correspondence/}\\
\texttt{letters/OO0766.pdf}}

\bibitem[13]{fantechi}
\textsc{Barbara Fantechi, Lothar G\"{o}ttsche, Luc Illusie, Steven L. Kleiman, Nitin Nitsure, Angelo Vistoli}, \rm{Fundamental algebraic geometry, Grothendieck's fga explained, Mathematical Surveys and Monographs vol. \textbf{123}, American Mathematical Society, Providence, RI, 2005.}

\bibitem[14]{GKZ}
\textsc{Israel M. Gelfand, Mikhail Kapranov, Andrei Zelevinsky}, \rm{Discriminants, resultants and multidimensional determinants, Birkh\"{a}user--Verlag, Boston, MA, 1994.} 

\bibitem[15]{gejapaper}
\textsc{Wolf-Dieter Geyer, Moshe Jarden}, Bounded realization of $l$-groups over global fields, 
\emph{Nagoya Math. J. \textbf{150} (1998) 13--62.}

\bibitem[16]{goldbach}
\textsc{Christian Goldbach}, \rm{Letter to Leonhard Euler dated 7th June 1742, 
Adolf P. Juskevic, Eduard Winter, eds., Leonhard Euler und Christian Goldbach, Briefwechsel 1729--1764, Akademie--Verlag, Berlin, 1965,
also available as Lettre XLIII from \texttt{http://www.math.dartmouth.edu/$\sim$euler/correspondence/}\\
\texttt{letters/OO0765.pdf}}

\bibitem[17]{AG-EGA4}
\textsc{Alexandre Grothendieck}, \rm{\'{E}l\'{e}ments de g\'{e}om\'{e}trique alg\'{e}brique (r\'{e}dig\'{e}s avec la collaboration de Jean Dieudonn\'{e}): IV. \'{E}tude locale des sch\'{e}mas et des morphismes de sch\'{e}mas, seconde partie, Publications math\'{e}matiques de l'I.H.\'{E}.S, tome \textbf{24}  (1965), 5--231.} 

\bibitem[18]{hayes}
\textsc{David R. Hayes}, \rm{A polynomial analog of the Goldbach conjecture, Bull. Amer. Math. Soc. \textbf{69} (1963), 115--116, Correction ibid. 493.}

\bibitem[19]{homma}
\textsc{Masaaki Homma}, \rm{Funny plane curves in characteristic $p>0$, Comm. Algebra \textbf{15} no. 7 (1987), 1469--1501.}

\bibitem[20]{hurwitz}
\textsc{Adolf Hurwitz}, \rm{\"{U}ber Riemann'sche Fl\"{a}chen mit gegebenen Verzweigungspunkten, Math. Ann. \textbf{39} (1891), 1--60 and Math. Werke, Band 1/XXI, Birkh\"{a}user, Basel, 1932.}

\bibitem[21]{Jou}
\textsc{Jean-Pierre Jouanolou}, \rm{Le formalisme du r\'{e}sultant, Adv. in Math. \textbf{90} no. 2 (1991), 117--263.}

\bibitem[22]{Lorenzini}
\textsc{Dino Lorenzini}, \rm{An invitation to arithmetic geometry, American Mathematical Society, Providence, RI, 1996.}

\bibitem[23]{nathanson}
\textsc{Melvyn B. Nathanson}, \rm{Additive number theory: the classical bases, Graduate Texts in Mathematics \textbf{164}, Springer--Verlag, New York, NY, 1996.}

\bibitem[24]{pollack-thesis}
\textsc{Paul Pollack}, \rm{Prime polynomials over finite fields, Ph.D. thesis, Dartmouth College, 2008.}

\bibitem[25]{pollack-article}
\textsc{Paul Pollack}, \rm{The exceptional set in the polynomial Goldbach problem, Int. J. Number Theory \textbf{7} no. 3 (2011), 579--591.}

\bibitem[26]{schmidt}
\textsc{Wolfgang M. Schmidt}, \rm{Equations over finite fields: an elementary approach, Lecture Notes in Mathematics \textbf{536}, Springer--Verlag, Heidelberg, 1976.}

\bibitem[27]{serre}
\textsc{Jean-Pierre Serre}, \rm{Topics in Galois Theory, Jones and Bartlett Publishers, Boston, MA, 1992.}

\bibitem[28]{vin}
\textsc{Ivan M. Vinogradov}, \rm{Representation of an odd number as a sum of three primes, Dokl. Akad. Nauk SSSR \textbf{15} no. 6--7 (1937), 291--294.}

\end{thebibliography}
\end{document}